\documentclass{amsart}
\usepackage{pb-diagram}
\newtheorem{lem}{Lemma}
\newtheorem{thm}{Theorem}

\newtheorem{cor}{Corollary}

\begin{document}

\title{Viewing AF-algebras as graph algebras}

\author{D. Drinen}
\subjclass{46L05, secondary 22A22}
\address{Department of Mathematics\\Arizona State University\\Tempe, AZ
85287}
\email{Drinen@asu.edu}

\begin{abstract}

Every AF-algebra $A$ arises as the $C^*$-algebra of a locally finite
pointed directed graph in the sense of Kumjian, Pask, Raeburn, and
Renault.  For AF-algebras, the diagonal subalgebra defined by
Str\v{a}til\v{a} and Voiculescu is consistent with Kumjian's notion of
diagonal, and the groupoid arising from a well-chosen Bratteli diagram for
$A$ coincides with Kumjian's twist groupoid constructed from a diagonal of
$A$.

\end{abstract}

\maketitle

\section{Introduction}

Constructing a Bratteli diagram (\cite{brat}) for an AF-algebra is a
well-known way of
associating a (directed) graph with a $C^*$-algebra.  Recently, Kumjian,
Pask, Raeburn, and Renault have been doing the opposite.  That is, they
introduced a notion of constructing a $C^*$-algebra from a directed graph.
For certain kinds of graphs, their strategy for doing so was to build from
the graph a locally compact $r$-discrete groupoid, and then use the theory
of Renault \cite{ren} to construct a $C^*$-algebra.  In this paper, we
show that these two procedures are compatible.

Suppose $E$ is a row-finite directed graph with no sinks.  In
\cite{kprr}, Kumjian, Pask, Raeburn, and Renault construct a locally
compact $r$-discrete groupoid $\mathcal{G}_E$ from $E$, and
study the associated $C^*$-algebra $C^*(\mathcal{G}_E)$. The class of
$C^*$-algebras that arise in
this way includes, up to Morita equivalence, Cuntz-Krieger algebras and 
AF-algebras.  In \cite{kpr}, Kumjian, Pask, and Raeburn define the
$C^*$-algebra $C^*(E)$ of any
row-finite directed graph $E$, and note that $C^*(E)$ coincides with
$C^*(\mathcal{G}_E)$ in case $E$ has no sinks.  They then show that
$C^*(E)$ is an AF-algebra if and only if $E$ has no (directed) loops.  The
main
result of this paper is that every AF-algebra arises in this way, at
least if we are willing to point the graph. That
is, every AF-algebra $A$ is the $C^*$-algebra of a row-finite
(pointed) directed
graph.  Not surprisingly, that directed directed graph will be a Bratteli
diagram for $A$.  An arbitrary Bratteli diagram for $A$, however, may not
produce the desired result, so we start in Section~\ref{sec-prelim} by
giving a procedure for building appropriate Bratteli diagrams.
Section~\ref{sec-mr} recalls the basics from \cite{kprr} and then proves
the main result.  

In Section~\ref{sec-cwd}, we explore the notion of the diagonal
subalgebra.  This is defined in \cite{s-v} for unital AF-algebras and in
\cite{diag} for general $C^*$-algebras.  We show that every AF-algebra
(unital or not)
contains a subalgebra which is diagonal in the sense of \cite{diag} and
that, not surprisingly, the two notions of diagonal coincide in the AF
case.  Further, our
construction of a groupoid from a  suitable Bratteli diagram for $A$
is compatible with Kumjian's construction in \cite{diag} of a
groupoid from a diagonal of $A$.

This paper was written while the author was a student at Arizona State
University, and it formed the author's masters thesis.  I would like
to
take this opportunity to express my gratitude to my
advisor, John Quigg, for all the patient help and guidance he provided.
I would also like to thank Alex Kumjian, David Pask, and Jack Spielberg  
for their helpful
discussions about the ideas in this paper.


\section{Preliminaries}


\label{sec-prelim}
Our conventions for AF-algebras come primarily from \cite{effros}.  Given
a
sequence of
finite-dimensional $C^*$-algebras $\{A_n\},$ and
injective $^*$- homomorphisms $\phi_n:A_n\hookrightarrow A_{n+1},$ 
we define $\lim_\rightarrow (A_n,
\phi_n)$ to be the $C^*$-algebraic direct limit of the inductive system 
\[
\begin{diagram}
\node{A_1}
	\arrow{e,t}{\phi_1}
\node{A_2}
	\arrow{e,t}{\phi_2}
\node{A_3}
	\arrow{e,t}{\phi_3}
\node{\cdots}
\end{diagram}
\]
The $C^*$-algebra $\lim_\rightarrow (A_n,\phi_n)$ is called an {\it
AF-algebra}.

Given an AF-algebra $$A=\lim_{\longrightarrow} (A_n,\phi_n),$$
let $p_n$ be the number of simple summands in the algebra $A_n$, and let
$[n,1]$, $[n,2]$,$\ldots$,$[n,p_n]$ denote the sizes of the corresponding
matrix
algebras.  That is, for each $n\in {\bf N},
A_n = \bigoplus_{i=1}^{p_n} A^n_i$, with $A^n_i \cong M_{[n,i]}$.
We will denote by $[K^{(n)}]$ the $p_{n+1}\times p_n$ matrix of
multiplicities corresponding to the map $\phi_n:A_n\rightarrow A_{n+1}$.
That is, the multiplicity of the embedding $\phi^n_{i,j}: A^n_j
\longrightarrow A^{n+1}_i$ is $k^{(n)}_{i,j}$.  For $m < n,$ define
$\phi_{n,m}=\phi_{n-1} \circ \phi_{n-2} \circ \cdots
\circ
\phi_{m}:A_m \longrightarrow A_n$.

Given this set-up, we construct a Bratteli diagram $G$ for $A$ as follows:
the vertex set $V$ for $G$ is defined by $V=\bigcup_{i=1}^{\infty} V^n$
and
$V^n=\{v_1^n, v_2^n, \ldots , v_{p_n}^n\}$.  For every $n\in N$, and for
every $i\leq p_{n+1}$ and $j\leq p_n$, put $k^{(n)}_{i,j}$ edges from
$v^n_j$ to $v^{n+1}_i$.  $G$ is now a row-finite (locally finite, in fact)
directed graph in the sense of \cite{kprr}.  There are, of course, many
Bratteli diagrams associated with a particular AF-algebra, but we will
call
the one constructed above the {\it Bratteli diagram corresponding to
the system $(A_n, \phi_n)$}. One more
bit of notation:  for $n>1, i\leq p_n$, put
$$\sigma^n_i = [n,i] - \sum_{j=1}^{p_{n-1}} k^{(n-1)}_{i,j} [n-1,j],$$
and put $\sigma^1_i=[1,i]$ for $i\leq p_1$.  
Think of $\sigma^n_i$ as the ``slack'' at vertex $i$ of level $n,$ the
extent
to which the embeddings of the algebras at level $n-1$ fail to fill up the 
algebra $A^n_i$.  Of course
$\sigma^n_i$ is a nonnegative integer for all $n,i$;  it will also be
convenient in what follows to arrange things so that $\sigma^n_i\leq
1$. The idea of the
following lemma is to replace a Bratteli diagram which has vertices with
too much slack with one that corresponds to the
same AF-algebra, but which has limited slack at all of its vertices.

\begin{lem}   
Every AF-algebra $A=\lim_\rightarrow (A_n, \phi_n)$ is isomorphic to
an AF-algebra $D=\lim_\rightarrow (D_n, \psi_n)$ which satisfies
$D_1={\bf C}$
and 
\begin{equation}
    \sigma^n_i\leq 1\,\hbox{for all } n>1\, \hbox{and all }i\leq
p_n.
\label{eq:slack}
\end{equation}
\end{lem}

Before we proceed with the proof, we present an example to guide the
reader.  Suppose we are confronted with the following undesirable Bratteli
diagram:

\[
\begin{diagram}
\node{3}
        \arrow{e}
\node{11}
        \arrow{e}
        \arrow{se}
\node{21}\\
\node{5} 
        \arrow{ne}
        \arrow{e,b}{2}
\node{10}
        \arrow{ne}
        \arrow{e}
\node{21}
\end{diagram}
\]

\noindent
First we fix the beginning:
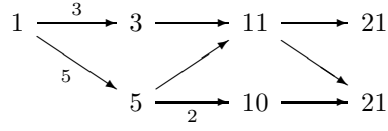
\[
\begin{diagram}
\node{1} 
        \arrow{e,t}{3}
        \arrow{se,b}{5}
\node{3} 
        \arrow{e}
\node{11}
        \arrow{e}
        \arrow{se}
\node{21}\\
\node{}
\node{5}
        \arrow{e,b}{2}
        \arrow{ne}
\node{10}
        \arrow{e}
\node{21}
\end{diagram}
\]

\noindent
Then we add a layer:
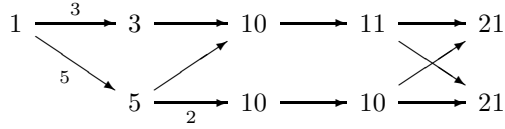
\[
\begin{diagram}
\node{1}
        \arrow{e,t}{3}
        \arrow{se,b}{5}
\node{3} 
        \arrow{e}
\node{10}
        \arrow{e}
\node{11}
        \arrow{e}
        \arrow{se}
\node{21}\\
\node{}
\node{5}
        \arrow{ne}
        \arrow{e,b}{2}
\node{10}
        \arrow{e}
\node{10}
        \arrow{ne}
        \arrow{e}
\node{21}
\end{diagram}
\]

\noindent
And another:

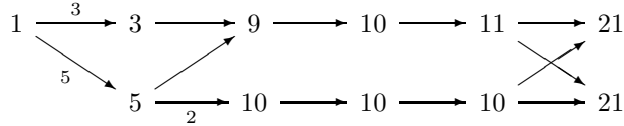
\[
\begin{diagram}   
\node{1}
        \arrow{e,t}{3}
        \arrow{se,b}{5}
\node{3}
        \arrow{e}
\node{9} 
        \arrow{e}
\node{10}
        \arrow{e}
\node{11}
        \arrow{e}
        \arrow{se}
\node{21}\\
\node{} 
\node{5}
        \arrow{ne}
        \arrow{e,b}{2}
\node{10}
        \arrow{e}
\node{10}
        \arrow{e}
\node{10}
        \arrow{ne}
        \arrow{e}
\node{21}
\end{diagram}
\]

\noindent
Further, we will fix the maps along the way so that, for
instance, the map $M_5 \hookrightarrow M_{11}$ in the
original
diagram
is the same as the composition of the three maps $M_5
\hookrightarrow M_9 \hookrightarrow M_{10} \hookrightarrow M_{11}$ in the
final diagram.

\begin{proof}  
Given $A$ with $A_1\neq {\bf C},$ define a new inductive system
$(B_n,\rho_n)$ as follows.  Let $B_1= {\bf C},$ and let $B_n=A_{n-1}$ for
$n>1.$  Define $\rho_n=\phi_{n-1}$ for $n>1,$ and let $\rho_1$ be the
unique unital homomorphism from ${\bf C}$ to $B_2$.  If $B$ is defined to
be $\lim_\rightarrow (B_n,\rho_n),$
then it is clear that $B\cong A.$  So we may assume, without loss of
generality, that $A_1={\bf C}.$  We now proceed to modify $A$ so that
it satisfies Condition (\ref{eq:slack}).

First let $n$ be the smallest integer
for which there is an $i\leq p_n$ with $\sigma^n_i > 1$.  Define $D_j =
A_j$ for $j=1\ldots n-1$ and let $\psi_j=\phi_j$ for $j=1\ldots n-2$.  Now
define $D_n$ as follows.  For $k=1 \ldots p_n$, if $\sigma^n_k \leq 1$,
put $D^n_k = A^n_k$.  If $\sigma^n_k > 1$, choose $D^n_k$ such
that $\bigcup_{j=1}^{p_{n-1}} \hbox{range}(\phi^{n-1}_{k,j}) \subset
D^n_k
\subset A^n_k$ and $D^n_k \cong M_{[n,k]-1}$.  Put $D_n =
\bigoplus_{k=1}^{p_n} D^n_k$.

Let $\psi_{n-1} = \phi_{n-1}$, but viewed as a map from $A_{n-1}$ to
$D_n$.  Now let $D_{n+1}=A_n$ and let $\psi_n :D_n
\hookrightarrow D_{n+1}$ denote the inclusion map.  It is clear that 
$\phi_{n-1}=\psi_n \circ \psi_{n-1}.$  

Now define $D_{n+j}=A_{n+j-1}$ for $j\geq 2$ and $\psi_{n+j} =
\phi_{n+j-1}$ for $j\geq 1$.  This gives a new inductive system, and thus
a new AF-algebra.  Now repeat this process as many times as necessary
(possibly countably many) and denote the resulting system also by $D = 
\lim_\rightarrow(D_n, \psi_n)$.  

It is clear that $D$ satisfies the
required condition, so it remains to show that $D\cong A.$  To that end,
we note that $D =
\lim_{\rightarrow} (D_{n_k}, \psi_{n_{k+1}, n_k})$ for any subsequence
$\{n_k\}$.  By construction, $A_1 = D_1$.  Inductively, for $k > 1$, we
can choose an $n_k > n_{k-1}$ such that $A_k = D_{n_k}.$  Because of the
way the $\psi_n$'s were constructed, the sequence $(D_{n_k},
\psi_{n_{k+1}, n_k})$ coincides with the sequence $(A_k, \phi_k)$.  Thus
$A \cong D$.

\end{proof}


\section{The Main Result}

\label{sec-mr}

Before we state the main result, we collect some notation and basic facts
from \cite{kprr}.  Given a directed graph $E,$ let $E^0$ denote the set of
vertices and $E^1$ the set of edges.  For $e \in E^1,$ let $s(e)$ and
$r(e)$ denote the source and range of $e,$ so $s,r: E^1 \longrightarrow
E^0.$  For $n \in {\bf N},$ let $E^n = \{e_1e_2 \ldots e_n\,|\,e_i
\in E^1, s(e_i) =
r(e_{i-1}), i=2 \ldots n \},$ the set of all paths in $E$ of length $n.$
Let $E^* = \bigcup_{n \in {\bf N}} E^n,$ the set of all finite paths, and
$E^\infty = \{e_1e_2\ldots\,|\,s(e_i)=r(e_{i-1})\},$ the infinite path
space of $E.$

For $x,y \in E^{\infty},\,k \in {\bf Z},$ say $x \sim_k y$ if and only if
$x_i = y_{i+k}$ for large $i.$  Define $\mathcal{G}_E,$ the path groupoid
of $E,$ by $\mathcal{G}_E = \{(x,k,y)\,|\,x \sim_k y\}.$  The groupoid
operations in $\mathcal{G}_E$ are as follows:
\begin{itemize}
\item $(x,k,y)^{-1} = (y,-k,x);$
\item $(x,k,y) \cdot (y,l,z) = (x,k+l,z).$
\end{itemize}

For $\alpha, \beta \in E^*$ with $r(\alpha) = r(\beta),$ define $Z(\alpha,
\beta) = \{(x,k,y)\,|\,x = \alpha x', y = \beta y'\}.$  The $Z(\alpha,
\beta)$'s are compact open $G$-sets in the topology they
generate.  Note that $r(x,k,y) = (x,0,x),$ so we identify $E^\infty$ with 
$\mathcal{G}_E^0,$ the unit space of  $\mathcal{G}_E.$  The relative
topology on $E^\infty$ as a subset of $\mathcal{G}_E$ agrees with the
product topology. 

Since the  $Z(\alpha, \beta)$'s are compact and open and generate the
topology of $\mathcal{G}_E,$ we have that $C^*(\mathcal{G}_E)$ is
generated by the characteristic functions of these sets.  That is,
$\hbox{span}\{1_{Z(\alpha, \beta)}\,|\,\alpha, \beta \in E^*, r(\alpha) = 
r(\beta)\}$ is dense in $C^*(\mathcal{G}_E).$  Note how two of these
characteristic functions multiply in $C^*(\mathcal{G}_E):$
 
$$1_{Z(\alpha , \beta )}*1_{Z(\gamma , \delta )} =
\left\{\begin{array}{cc}
        1_{Z(\alpha \epsilon , \delta )} & \hbox{ if } \gamma = \beta
\epsilon\\
        1_{Z(\alpha , \delta \epsilon )} & \hbox{ if } \beta = \gamma 
\epsilon\\ 
	0                               & \hbox{ otherwise,}
       \end{array}
\right. $$

\noindent
that is, the product is $0$ unless one of $\beta, \gamma$ continues the
other.

One more piece of notation:  suppose $S$ is any subset of $E^0$.
Considering only those paths whose source is in $S$ determines an open
subset of $\mathcal{G}_E^0$.  We will refer to the reduction of the
groupoid
$\mathcal{G}_E$ to this open set as $\mathcal{G}_E|_S$, slightly abusing
the notation.  Distinguishing a
set of vertices of the graph is known as {\it pointing} the graph.  In
what follows, we will be dealing with path groupoids reduced to such
subsets of the unit space.  Because $S$ is open, the Haar system for
$\mathcal{G}_E$ restricts to a Haar system for $\mathcal{G}_E|_S.$  Also,
$\{Z(\alpha, \beta)\,|\,s(\alpha), s(\beta) \in S\}$
is a base of compact open $G$-sets for the topology of $\mathcal{G}_E|_S,$
so the span of the characteristic functions of those sets is dense in
$C^*(\mathcal{G}_E|_S).$

We now state the main result:

\begin{thm}
\label{main-thm}
Given an AF-algebra $$A=\lim_{\longrightarrow} (A_n,\phi_n),$$ there is a
Bratteli diagram $E$ for $A$ and a set $S$ of vertices such that
\hbox{$C^*(\mathcal{G}_E|_S)\cong A$}. 
\end{thm}

\begin{proof}  First use the above procedure (if necessary) to construct
a new inductive system for the AF-algebra such that the Bratteli
diagram corresponding to that inductive system satisfies Condition
(\ref{eq:slack}).  That is, build a Bratteli diagram for $A$ which has no
excess slack.
Denote by $E$ the Bratteli diagram corresponding to this
inductive system.
Next,  define $S$ to be the set of all vertices where the slack is 1, 
along with the vertex at the top level, and for each $n$ and $i$ define
$F_i^n$ to be the set of all finite paths in $E$ that start in $S$ and end
at the vertex
$v^n_i$.  That is, $S:=\{v^n_i\,|\,\sigma_i^n=1\}\cup
\{v_1^1\}$ and $F^n_i:=\{\alpha \in
F(E)\,|\,s(\alpha ) \in S, r(\alpha )=v^n_i\}.$  Let $|F_i^n|$ denote the
cardinality of $F_i^n.$  

Now we claim that for
every integer $n,$
\begin{equation}
|F^n_i|=[n,i] \hbox{ for every } i\leq p_n.  
\label{eq:numpaths}
\end{equation}
The proof is
by induction.  First note that the case $n=1$ is taken care of by the
construction of $A.$  That is, since $p_1=1$ and $[1,1]=1,$ the set
$F_1^1$ has only one element, the zero-length path $v_1^1.$  Now suppose
that (\ref{eq:numpaths}) is
true for $n=l.$  Fix $i \leq p_{l+1}.$  If
$\sigma^{l+1}_i=0,$
then $v^{l+1}_i$ is not an element of $S,$ so every path starting in $S$
and
ending at $v^{l+1}_i$ must have come through some vertex in level $l.$
There are $p_l$ such vertices, the $j$th of which, by the induction
hypothesis, being the range of $[l,j]$ paths emanating from $S.$  There
are
$k^{(l)}_{i,j}$ different ways to get from $v^l_j$ to $v^{l+1}_i$, so we
must have 
$$|F^{l+1}_i| = \sum_{j=1}^{p_l} k^{(l)}_{i,j} [l,j] = [l+1,i] -
\sigma^{l+1}_i = [l+1,i],$$
as desired.  If $\sigma^{l+1}_i=1,$ then the calculations are similar:  we
end up with $[l+1,i]-1$ paths which pass through the preceding
level.  Adding in the zero-length path $v^{l+1}_i,$ which is in $S,$
brings us to the right number.  This establishes the claim.

Now define, for any positive integer $n$ and any $i \leq p_n,\,
B^n_i:=\hbox{span}\,\{1_{Z(\alpha , \beta)}\,|\,\alpha , \beta \in
F^n_i\} \subset C^*(\mathcal{G}_E|_S).$  Note that, given two generators
of $B^n_i,$ say $1_{Z(\alpha , \beta )}$ and $1_{Z(\gamma , \delta )},$
because $\beta$ and $\gamma$ end at the same vertex, one cannot continue the
other unless they coincide, so we have:
$$1_{Z(\alpha , \beta )}*1_{Z(\gamma , \delta )} = 
\left\{\begin{array}{cc}
	1_{Z(\alpha , \delta )}	& \hbox{ if } \beta =\gamma \\
        0				& \hbox{ otherwise. }
       \end{array}
\right. $$
\noindent	    
Also, $1_{Z(\alpha, \beta)}^* = 1_{Z(\beta, \alpha)}.$  So the elements of
$B^n_i$ behave
like matrix units, and there are, by the above claim, $[n,i]^2$ of them,
so $B^n_i\cong M_{[n,i]}.$  Further, if $i\neq j,$ $B^n_iB^n_j=0$
because if a path ends at $v^n_i,$ it cannot be a continuation of (or be
continued by) a path ending at $v^n_j.$  Thus, $B_n:= 
\hbox{span}B^n_i$ coincides with the direct sum $\bigoplus_{i=1}^{p_n}
B^n_i,$
and we have $$B_n \cong \bigoplus_{i=1}^{p_n} M_{[n,i]} \cong A_n.$$

Note now that each of these algebras embeds naturally into the next; for
instance, 
$$1_{Z(\alpha , \beta )}=\sum_{\{e \in E^1\,|\,s(e)=r(\alpha)\}}
1_{Z(\alpha
e,\beta e)}$$ gives the
inclusion $\iota_n$ of $B_n$ into $B_{n+1}$.  The multiplicity of
the embedding of $B^n_j$ into $B^{n+1}_i$ is the number of edges from
$v^n_j$ to $v^{n+1}_i,$ which is $k^{(n)}_{i,j}.$  All that remains is to
show that isomorphisms can be chosen between $B_n$ and $A_n$ such that the
following diagram commutes: 

\[
\begin{diagram}
\node{A_1={\bf C}}
        \arrow{e,t}{\phi_1}
\node{A_{2}}
	\arrow{e,t}{\phi_2}
\node{\cdots}
	\arrow{e}
\node{A_n}
	\arrow{e,t}{\phi_n}
\node{A_{n+1}}
	\arrow{e}
\node{\cdots}\\
\node{B_1={\bf C}}
        \arrow{e,b}{\iota_1}
	\arrow{n,l}{=}
\node{B_2}
        \arrow{e,b}{\iota_2}
	\arrow{n,l}{\cong}
\node{\cdots}
        \arrow{e}
\node{B_n}
        \arrow{e,b}{\iota_n}
	\arrow{n,l}{\cong}
\node{B_{n+1}}
        \arrow{e}
	\arrow{n,l}{\cong}
\node{\cdots}
\end{diagram}
\]

First, for each $n,$ choose any isomorphism $\hat{\rho}_n$ from $B_n$
to $A_n.$  Since $\phi_1$ has the same multiplicities as $\iota_1 \circ
\hat{\rho}_2,$ there is a unitary $u_2 \in A_2$ such that $\phi_1 =
\hbox{Ad}u_2
\circ \hat{\rho}_2 \circ \iota_1.$  So let the isomorphism $\rho_2$
between
$B_2$ and
$A_2$ be $\hbox{Ad}u_2 \circ \hat{\rho}_2.$  Inductively, if isomorphisms
$\{\rho_j\}_{j=2}^n$ have been chosen to make the first $n-1$ squares of
the diagram commute, choose
a unitary $u_{n+1} \in A_{n+1}$ such that $\hbox{Ad}u_{n+1} \circ
\hat{\rho}_{n+1} \circ \iota_n = \phi_n \circ \rho_n,$ as we may, since 
the maps $\phi_n \circ \rho_n$ and $\hat{\rho}_{n+1} \circ \iota_n$
have the same multiplicities, and let $\rho_{n+1} = \hbox{Ad}u_{n+1}
\circ \hat{\rho}_{n+1}.$ 

So if $B$ is defined to be $\lim_\rightarrow (B_n,\iota_n),$
we have $$B\cong \lim_{\longrightarrow} (A_n,\phi_n)=A.$$
Since $\bigcup_{n=1}^\infty B_n = \hbox{span}\{1_{Z(\alpha
,\beta)}\,|\,\alpha,\beta\in 
F_i^n \hbox{ for some $n,i$}\,\}$ generates
$C^*(\mathcal{G}_E|_S),$ it follows that $A\cong
C^*(\mathcal{G}_E|_S).$  
\end{proof}

Note that it is necessary to correctly modify and correctly point the
Bratteli diagram, as the following simple example shows:  $M_3$ can be
thought of as an AF-algebra with the following Bratteli diagram:
 
\[
\begin{diagram}
\node{1}
	\arrow{e}
\node{2}
	\arrow{e}
\node{3}
	\arrow{e}
\node{3}
	\arrow{e}
\node{\cdots}
\end{diagram}
\]

\noindent
If we turn this into a directed graph in the obvious way and fail to point
it, we obtain the following graph:

\[
\begin{diagram}
\node{\cdot}
        \arrow{e}
\node{\cdot}   
        \arrow{e}
\node{\cdot}
        \arrow{e}
\node{\cdot}
        \arrow{e}
\node{\cdots}
\end{diagram}
\]

\noindent
of which the $C^*-$algebra is $\mathcal{K},$ the compact operators.  If we
point the graph, for instance, at only one vertex, we obtain the graph:

\[
\begin{diagram}
\node{*}
        \arrow{e}
\node{\cdot}   
        \arrow{e}
\node{\cdot}
        \arrow{e}
\node{\cdot}
        \arrow{e}
\node{\cdots}
\end{diagram}
\]

\noindent
whose $C^*-$algebra is {\bf C}.  Applying the procedure of Section 2 to
the Bratteli diagram in question would produce the following pointed
graph:

\[
\begin{diagram}  
\node{*}
        \arrow{e}
\node{*}
        \arrow{e}
\node{*}
        \arrow{e}
\node{\cdot}
        \arrow{e}
\node{\cdots}
\end{diagram}
\]

\noindent
which has the desired algebra, $M_3,$ as a $C^*-$algebra.

Note that there are other procedures for modifying and pointing the
Bratteli diagram that would yield the same result.  One such (communicated
to us by Alex Kumjian) would be to
add a vertex to every level below which there is a vertex with slack
greater than 0, then add the appropriate number of edges from the new
vertex to the ones with slack.  That is, if $\sigma_{n,i} > 0$ for some
$i \leq p_n$, add a
vertex to level $n-1$ and add $\sigma_{n,i}$ edges from the new vertex to
the vertex $v^n_i$.  Then point this Bratteli diagram at the added
vertices.  This approach is somewhat more natural because every pointed
vertex is a source, but it makes for unnatural Bratteli
diagrams.

Also note that, in case $A$ is unital, it is possible to construct a
Bratteli diagram for $A$ which has no slack at all and which has a single
vertex
at the top level.  In this case, it is only necessary to point the top
vertex.

%
%

\section{Connection with the Diagonal}

%
%
\label{sec-cwd}

In \cite{diag}, Kumjian defines the notion of a diagonal subalgebra of a
$C^*$-algebra $A$ as follows.  An abelian subalgebra $B$ of $A$ is said to
be
{\it diagonal in $A$} if it contains a positive element which is strictly
positive in $A$ and if there exists a faithful conditional expectation
$P:A \longrightarrow B$ such that $N_f(B)$ spans a dense subset of
$\hbox{ker} P.$  Here, $N_f(B)$ denotes the {\it free normalizers} of $B$
in
$A$, namely the set of all $a \in A$ such that
$a^*Ba \cup  aBa^* \subset B$ and $a^2=0$.  He defines a {\it twist} as a
proper ${\bf T}$-groupoid $\Gamma$ such that $\Gamma / {\bf T}$ is an
$r$-discrete equivalence relation, and shows that there is a one-to-one
correspondence
between twists and diagonal pairs of $C^*$-algebras.  We briefly recall
the basics here.

Suppose $\Gamma$ is a twist.  Define $$E(\Gamma) = \{f \in
C_c(\Gamma)\,|\, f(t\gamma) = tf(\gamma) \hbox{ for all }t \in {\bf T},
\gamma \in \Gamma\}$$  

\noindent and

$$D(\Gamma) = \{f \in E(\Gamma)\,|\,\hbox{supp}\,f \subset
{\bf T}\Gamma^0\},$$

\noindent where ${\bf T}\Gamma^0$ denotes the istotropy group bundle of
$\Gamma$.  

$E(\Gamma)$ is then a $^*$-algebra with a
distinguished abelian subalgebra $D(\Gamma)$.  Further, $D(\Gamma) \cong
C_c(\Gamma^0)$.  $E(\Gamma)$ becomes a pre-Hilbert $D(\Gamma)$-module; the
completion, which Kumjian denotes by $\mathcal{H}(\Gamma)$, is a Hilbert
$C_0(\Gamma^0)$-module.

Kumjian then
constructs a $^*$-homomorphism $\pi : E(\Gamma) \longrightarrow
L(\mathcal{H}(\Gamma))$ such that $\pi(f)g = fg$ (the convolution
product) for all $f,g \in
E(\Gamma)$ and
defines $A(\Gamma)$ to be the closure of
$\pi (E(\Gamma))$ and $B(\Gamma)$ to be the closure of $\pi
(D(\Gamma))$.  It turns
out that $B(\Gamma)$ is diagonal in $A(\Gamma)$, and this shows that
every twist gives rise to a diagonal pair of $C^*$-algebras.

Conversely, given $C^*$-algebras $A$ and $B$ with $B$ diagonal in $A$,
Kumjian constructs a twist $\Gamma$ such that $A \cong A(\Gamma)$ and $B
\cong B(\Gamma)$, giving a bijective correspondence between
twists and diagonal pairs.

The following theorem is perhaps known to experts, but we could not find
it in the literature.

\begin{thm}
\label{diagthm}
Let $\mathcal{G}$ be a Hausdorff, amenable, $r$-discrete equivalence
relation.  Then $(C^*(\mathcal{G}),
C_0(\mathcal{G}^0))$ is a diagonal pair in Kumjian's sense, and
$\mathcal{G} \cong \Gamma / {\bf T}$, where $\Gamma$ is Kumjian's
associated twist.
\end{thm}

\begin{proof}
By amenability, principality, and \cite[Proposition II.4.7(ii)]{ren},
$C_0(\mathcal{G}^0)$ is a maximal abelian subalgebra of
$C^*(\mathcal{G})$.  By amenability and \cite[Proposition II.4.8]{ren},
the map $P:
C_c(\mathcal{G}) \longrightarrow C_0(\mathcal{G}^0)$ given by restriction
to the unit space extends to a faithful conditional expectation (also
denoted by $P$) of $C^*(\mathcal{G})$ onto $C_0(\mathcal{G}^0)$.  To
see that $C_0(\mathcal{G})$ is diagonal in $C^*(\mathcal{G})$, we must
show that $\hbox{Ker}P = \overline{\hbox{span}}N_f(C_0(\mathcal{G}^0))$.
To
do so, we adapt the proof of \cite[Lemma 2.12]{diag}.
 
\begin{lem}
\label{denselem}
$\hbox{Ker}P = \overline{\hbox{Ker}P \cap C_c(\mathcal{G})}.$
\end{lem}

\begin{proof}
Fix $\epsilon > 0$ and $x \in \hbox{Ker}P$.  There is $f \in
C_c(\mathcal{G})$ such that $||f - x|| < \epsilon / 2$.  Thus, because $P$
is contractive, 
$||f|_{\mathcal{G}^0}|| < \epsilon / 2$.  So
letting $g = f - f|_{\mathcal{G}^0}$, and noting that $\mathcal{G}^0$ is
clopen (because $\mathcal{G}$ is $r$-discrete), we have $g \in
\hbox{Ker}P \cap C_c(\mathcal{G})$ and  $||g - x|| < \epsilon$.
\end{proof}

\begin{lem}
\label{disjointlem}
If $x \in \mathcal{G}$ and $x$ is not a unit, then there exists an open
$G$-set $S$ such that $x \in S$ and $s(S) \cap r(S) = \emptyset$
\end{lem}

\begin{proof}
Since $\mathcal{G}$ is Hausdorff and $s(x) \neq r(x)$, there are disjoint 
open subsets
$U_1$ and $U_2$ of the unit space such that $s(x) \in U_1$ and $r(x) \in
U_2$.  Also, $x$ must be contained in some open $G$-set $S_1$.
So $S := S_1 \cap \mathcal{G}_{U_1}^{U_2}$ is an open
$G$-set which contains $x$ and whose source is disjoint from its range. 
\end{proof}

By Lemma~\ref{denselem}, it
suffices to show that for every $f \in \hbox{Ker}P \cap
C_c(\mathcal{G})$, there exist $g_1$, $g_2$, $\ldots$, $g_n \in
N(C_0(\mathcal{G}^0))$ such that $f = \sum g_i$ and $g_i^2 = 0$ for
every $i$.  

Since $f \in \hbox{Ker}P$, $f|_{\mathcal{G}^0} = 0$.  Thus, by
Lemma~\ref{disjointlem} and the fact that $f$ has compact support, there
exist finitely many open $G$-sets $S_1$, $S_2$, $\ldots$, $S_n$
such that $\hbox{supp}(f) \subset \bigcup S_i$ and $s(S_i) \cap r(S_i) =
\emptyset$ for every $i$. 

Let $\{h_i\}$ be a partition of unity subordinate to
$\{S_i\}$, and for every $i$ define $g_i = fh_i$.  It is clear that $f =
\sum g_i$. We
claim that $g_i \in N_f(C_0(\mathcal{G}^0))$ for every $i$.  First, note
that $\hbox{supp}\,(g_i^2) \subset (\hbox{supp}\, g_i)^2 \subset S_i^2 =
\emptyset$, so $g_i^2 = 0$.  To show that each $g_i$ is a normalizer of
$C_0(\mathcal{G}^0)$, it
suffices to show that $g_i$ normalizes $C_c(\mathcal{G}^0)$.  To see this
note that, if $h \in C_c(\mathcal{G}^0)$, then $\hbox{supp}\,(g_i h g_i^*) 
\subset \hbox{supp}\,(g_i) \, \hbox{supp}\,(h) \, \hbox{supp}\,(g_i^*)
\subset S_i \mathcal{G}^0 S_i^{-1} = r(S_i)$.  This last equality comes
from the fact that $S_i$ is a $G$-set.  So we have $\hbox{supp}\,(g_i h
g_i^*) \subset \mathcal{G}^0$.  Clearly, $g_i h g_i^*$ is continous and
has compact support, so $g_i h g_i^* \in C_c(\mathcal{G}^0)$.  Similarly
for $g_i^* h g_i$, so $g_i$ is a free normalizer.

This establishes that $C_0(\mathcal{G}^0)$ is diagonal in
$C^*(\mathcal{G})$.  We now proceed to show that $\mathcal{G} \cong \Gamma
/{\bf T}$.

Given $\mathcal{G}$, we construct the trivial twist $\Gamma_1$ over
$\mathcal{G}$.  That is, let $\Gamma_1 = \mathcal{G} \times {\bf T}$.  The
${\bf T}$-groupoid operations on $\Gamma_1$ are as follows:

\begin{itemize}
\item $(g_1,t_1)(g_2,t_2) = (g_1g_2,t_1t_2);$
\item $(g,t)^{-1} = (g^{-1}, \overline{t});$
\item $s(g,t) = (g,t)^{-1}(g,t) = (s(g),1)$ (hence we identify the unit
space of $\Gamma_1$ with the unit space of $\mathcal{G}$);
\item $t_1(g,t_2) = (g,t_1t_2).$
\end{itemize}

Note that $\mathcal{G}$ can be viewed as a subgroupoid of $\Gamma_1$ via
$g \mapsto (g,1)$.  Thus, given a continuous equivariant function $f$
with compact
support
in $\Gamma_1$ (i.e. an element of Kumjian's $E(\Gamma_1)$), we can define
$\hat{f} \in C_c(\mathcal{G})$ by $\hat{f} = f|_{\mathcal{G}}$.  It is
easily checked that $\hat{ }:E(\Gamma_1) \longrightarrow C_c(\mathcal{G})$
is an isomorphism of $^*$-algebras.  

Thus, since $C_c(\mathcal{G})$ is dense in both $A(\Gamma_1)$ and
$C^*(\mathcal{G})$, in order to see that $A(\Gamma_1) \cong
C^*(\mathcal{G}),$ it suffices to show: 

\begin{equation}
\nonumber
||f||_{A(\Gamma_1)} 
= ||f||_{C^*(\mathcal{G})} \qquad \hbox{for every }f \in
C_c(\mathcal{G}).
\end{equation}

Because $\mathcal{G}$ is amenable, $||f||_{C^*(\mathcal{G})} = \sup_{u
\in \mathcal{G}^0}
||\hbox{Ind}\epsilon_u(f)||$, where $\hbox{Ind}\epsilon_u$ is the
representation of $C_c(\mathcal{G})$ on $L^2(\nu_u^{-1})$ given by
$\hbox{Ind}\epsilon_u(f)\xi(x) = \int f(y) \xi (y^{-1}x)
d\lambda^{r(x)}(y)$.  Here, $\nu_u$ denotes the measure on
$\mathcal{G}$
induced by the measure $\epsilon_u$ on $\mathcal{G}^0$, given by $\nu_u =
\int_{v \in \mathcal{G}^0} \lambda^v d\epsilon_u(v)$.  The measure
$\nu_u^{-1}$ is given by $\int \phi(x) d\nu_u^{-1}(x) = \int \phi(x^{-1})
d\nu_u(x).$  

Now, for $u \in \mathcal{G}^0$, consider the representation $\delta_u :
C_0(\mathcal{G}^0) \longrightarrow {\bf C}$ given by $\delta_u(f) = f(u)$.
Since $\mathcal{H}(\Gamma_1)$ is a (right) Hilbert
$C_0(\mathcal{G}^0)$-module with a left action of $A(\Gamma_1)$ by
adjointable operators, we can induce the
representation $\delta_u$ up to a representation $\hbox{Ind}\delta_u$ of
$A(\Gamma_1)$ on the Hilbert space $\mathcal{H}(\Gamma_1)
\otimes_{C_0(\mathcal{G}^0)} {\bf C}$.  We claim that for every $u \in
\mathcal{G}^0$, $\hbox{Ind}\delta_u$ is unitarily equivalent to
$\hbox{Ind}\epsilon_u$.

Define $v : \mathcal{H}(\Gamma_1) \odot {\bf C}
\longrightarrow L^2(\nu_u^{-1})$ by $v(\sum g_i \otimes \lambda_i) = \sum
\lambda_i g_i$.  We have

\begin{eqnarray}
\langle \sum_ig_i \otimes \lambda_i , \sum_j h_j \otimes \mu_j
\rangle_{\mathcal{H}(\Gamma_1) \otimes {\bf C}}
& = & \sum_{i,j}\langle g_i \otimes \lambda_i , h_j \otimes \mu_j\rangle
\nonumber \\
& = & \sum_{i,j}\overline{\lambda_i}\mu_j g_i^*h_j(u) \nonumber \\
& = & \sum_{x \in \mathcal{G}^u} \sum_{i,j} \overline{\lambda_i}\mu_j 
\overline{g_i(x^{-1})}h_j(x^{-1}) \nonumber \\
& = & \sum_{x \in \mathcal{G}^u} \overline{v(\sum_i g_i \otimes
\lambda_i)(x^{-1})} v(\sum_j h_j \otimes \mu_j)(x^{-1}) \nonumber \\
& = & \langle v(\sum_i g_i \otimes \lambda_i) , v(\sum_j h_j \otimes
\mu_j)\rangle_{L^2(\nu_u^{-1})} \nonumber
\end{eqnarray}

Thus $v$ can be extended to a unitary (also denoted by $v$) from
$\mathcal{H}(\Gamma_1) \otimes_{C_0(\mathcal{G}^0)} {\bf C}$ to
$L^2(\nu_u^{-1})$.  Also, for all $f \in C_c(\mathcal{G})$,

\begin{eqnarray}
v\,\hbox{Ind}\delta_u(f)(\sum g_i \otimes \lambda_i)
& = & v(\sum fg_i \otimes \lambda_i) \nonumber \\
& = & \sum \lambda_i fg_i \nonumber \\
& = & \sum \lambda_i \hbox{Ind}\epsilon_u(f)g_i \nonumber \\
& = & \hbox{Ind}\epsilon_u(f)\,v(\sum g_i \otimes \lambda_i). \nonumber
\end{eqnarray}

\noindent Thus $\hbox{Ind}\delta_u$ is unitarily equivalent to
$\hbox{Ind}\epsilon_u$ for every $u \in \mathcal{G}^0.$

Since $\{\delta_u\}$ is a separating family of representations of
$C_0(\mathcal{G}^0)$, $\{\hbox{Ind}\delta_u\}$ is a separating family of
representations of $A(\Gamma_1)$.  Thus, for $f \in C_c(\mathcal{G})$,

$$||f||_{A(\Gamma_1)} = \sup_{u \in \mathcal{G}^0}
||\hbox{Ind}\delta_u(f)|| = \sup_{u \in \mathcal{G}^0} 
||\hbox{Ind}\epsilon_u(f)|| = ||f||_{C^*(\mathcal{G})},$$ 

\noindent as desired.

Since the correspondence between diagonal pairs and twists is bijective
(up to isomorphism), we must have $\Gamma_1 \cong \Gamma$, where $\Gamma$
is the twist constructed from the diagonal pair $(C^*(\mathcal{G}),
C_0(\mathcal{G}^0))$.  Since $\mathcal{G} \cong \Gamma_1 / {\bf T}$, we
get $\mathcal{G} \cong \Gamma / {\bf T}$.
\end{proof}

\begin{cor}
Every AF-algebra contains a subalgebra which is diagonal in Kumjian's
sense.
\end{cor}

\begin{proof}
By Theorem~\ref{main-thm}, every AF-algebra $A$ is the $C^*$-algebra of a
Hausdorff, $r$-discrete equivalence relation, which is amenable by
\cite[Corollary 5.5]{kprr}.  Thus, by Theorem~\ref{diagthm}, $A$ contains
a diagonal subalgebra. 
\end{proof}

In \cite{s-v}, Str\v{a}til\v{a} and Voiculescu defined a notion of a
diagonal
subalgebra for AF-algebras, and it is not surprising that the two notions
of diagonal coincide in the AF case.  We briefly review the
Str\v{a}til\v{a}-Voiculescu set-up.  

Given a unital AF-algebra $A = \overline{\bigcup_{n=0}^{\infty}A_n}$,
where
it is assumed that $A_0 \cong {\bf C}$, Str\v{a}til\v{a} and Voiculescu
inductively define an ascending sequence of abelian $C^*$-subalgebras
$\{C_n\}$ as follows:

\begin{itemize}
\item $C_0 = A_0;$
\item given $C_n$, define $C_{n+1} = \hbox{span}\,C_nD_{n+1}$, where
$D_{n+1}$ is an arbitrary maximal abelian self-adjoint subalgebra of
$A_{n+1} \cap A_n^{'}.$
\end{itemize}

\noindent The diagonal subalgebra $C$ is defined to be the closure of the
union of the $C_n$'s.  In \cite{s-v}, only unital AF-algebras were
considered, but the definition of $C$ makes sense in the non-unital case,
and it is in that more general setting that the following theorem is
proved.

\begin{thm}
Let $A$ be an AF-algebra.  Then Str\v{a}til\v{a} and Voiculescu's
subalgebra
$C$ is diagonal in Kumjian's sense.
\end{thm}

\begin{proof}
In light of Theorem~\ref{main-thm}, we can view $A$ as
$C^*(\mathcal{G}_E|_S)$ for some Bratteli diagram $E$ for $A$ and some
subset $S$ of vertices of $E$.  We claim that the $D_n$'s can be chosen
such that for each $n$, $C_n = \hbox{span}\,\{1_{Z(\alpha,
\alpha)}\,|\,r(\alpha) \in V^n\}$ (recall that $V^n$ denotes the $n$th
level of vertices of the Bratteli diagram).

Suppose $C_0$, $C_1, \ldots ,$ $C_n$ are as in the above claim.  For each
edge $e$ with $r(e) \in V^{n+1}$, define $f_e = \sum 1_{Z(\alpha,
\alpha)}$, where the sum runs over all paths $\alpha$ whose last edge is
$e$ (we also make the implicit assumption, here and from now on, that all 
paths start in the distinguished set $S$).  Also, for every $v \in V^{n+1}
\cap S$, define $f_v = 1_{Z(v,v)}$.  Now define

$$D_{n+1} = \hbox{span}\,\{f_e, f_v\,|\,r(e) \in V^{n+1}, v \in V^{n+1}
\cap S \}.$$

\begin{lem}
$D_{n+1}$ is a MASA in $A_{n+1} \cap A_n^{'}$.
\end{lem}  

\begin{proof}
Clearly $D_{n+1}$ is abelian.  To show that it is a MASA in $A_{n+1} \cap
A_n^{'}$, we will show that $x \in A_{n+1} \cap A_n^{'} \cap D_{n+1}^{'}$
implies $x \in D_{n+1}.$  So fix $x \in A_{n+1} \cap A_n^{'} \cap
D_{n+1}^{'}$.  We can write $x = \sum \lambda_{\alpha, \beta}
1_{Z(\alpha, \beta)}$, where the sum runs over all $\alpha, \beta$ with
$r(\alpha) = r(\beta) \in V^{n+1}$ and each pair appears only once.

First suppose that there exist $\alpha_0 \neq \beta_0$ such that
$\lambda_{\alpha_0, \beta_0} \neq 0$.  Then one of the following three
cases
holds:

\noindent Case 1:  One of $\alpha_0$ or $\beta_0$ is a zero-length path.

\noindent Case 2:  $\alpha_0$ and $\beta_0$ have the same last edge (hence
differ somewhere before the last edge).

\noindent Case 3:  $\alpha_0$ and $\beta_0$ have different last edges.

\noindent For Case 1, assume (without loss of generality) that $|\alpha_0|
= 0$ and let $y = 1_{Z(\alpha_0, \alpha_0)} \in D_{n+1}$.  For Case 2,
write $\alpha_0 = \gamma_0e$ where $e$ is an edge, and let $y =
1_{Z(\gamma_0, \gamma_0)} \in A_n$.  For Case 3, denote by $e_0$ the last
edge of $\alpha_0$, and let $y = f_{e_0} \in D_{n+1}$.  

In all three cases, $x$ and $y$ should commute, but the reader may
tediously verify that $1_{Z(\beta_0,
\alpha_0)}xy1_{Z(\beta_0, \alpha_0)} = 0$ and $1_{Z(\beta_0,             
\alpha_0)}yx1_{Z(\beta_0, \alpha_0)} = \lambda_{\alpha_0, \beta_0}
1_{Z(\beta_0, \alpha_0)} \neq 0$.  Thus $x$ has no ``off-diagonal''
entries.
So we can write $x = \sum\lambda_\alpha 1_{Z(\alpha, \alpha)}$. 

Now suppose there exist $\alpha_1 \neq \alpha_2$ such that $\alpha_1,
\alpha_2$ have the same last edge, but $\lambda_{\alpha_1} \neq
\lambda_{\alpha_2}$.  For $i = 1,2$, write $\alpha_i = \gamma_ie$ where
$e$ is an edge, and set $y = 1_{Z(\gamma_1, \gamma_2)} \in A_n$.  The
reader may check that $1_{Z(\alpha_1, \alpha_1)}xy \neq 1_{Z(\alpha_1,
\alpha_1)}yx$.  Hence $xy \neq yx$, violating the assumption that $x \in
A_n^{'}$.  

Thus $x = \sum \lambda_\alpha 1_{Z(\alpha, \alpha)} \in A_{n+1} \cap
A_n^{'} \cap D_{n+1}^{'}$ implies that $\lambda_{\alpha_1} =
\lambda_{\alpha_2}$ whenever the last edges of $\alpha_1$ and $\alpha_2$
agree, which implies that $x \in D_{n+1}$.    
\end{proof}

Now, $C_{n+1} = \hbox{span}\,C_nD_{n+1} \subset
\hbox{span}\,\{1_{Z(\alpha, \alpha)}\,|\,r(\alpha \in V^{n+1}\}$.  To show
the reverse inclusion, fix $\alpha$ with $r(\alpha) \in V^{n+1}$.  If
$|\alpha| = 0$, $1_{Z(\alpha, \alpha)} \in D_{n+1} \subset C_{n+1}$.  If
not, write $\alpha = \gamma e$ where $e$ is an edge and note that
$1_{Z(\alpha, \alpha)} = f_e1_{Z(\gamma, \gamma)} \in
\hbox{span}\,\{C_n, D_{n+1}\}$.  Thus, with the $D_n$'s chosen in this
way, $C_n = \hbox{span}\,\{1_{Z(\alpha, \alpha)}\,|\,r(\alpha) \in
V^{n}\}$ for every $n$, so $C = \overline{\hbox{span}}\,\{1_{Z(\alpha,
\alpha)}\} =
C_0((\mathcal{G}_E|_S)^0)$, which is diagonal in Kumjian's sense.

Now we have shown that the $C_n$'s can be chosen in such a way that $C$ is
diagonal.  Now, suppose $\tilde{C}$ is another Str\v{a}til\v{a}-Voiculescu
diagonal.  Then there exists an automorphism $\alpha$ of $A$ which maps
$C$ to $\tilde{C}$.  In fact, this automorphism can be chosen to be
approximately inner (i.e. there exists a sequence $\{u_k\}$ of unitaries
in $A$ such that $\alpha(a) = \displaystyle{\lim_{k \rightarrow \infty}}
u_k a u_k^*$.  See \cite{power}.).  Defining $\tilde{P} = \alpha \circ P \circ \alpha^{-1}$, the reader may
check that $\tilde{P}$ is a faithful conditional expectation from $A$ onto
$\tilde{C}$, and that $\hbox{ker}\,\tilde{P} = \alpha(\hbox{ker}\,P)
= \overline{\hbox{span}}\,N_f(\tilde{C})$.  Thus $\tilde{C}$ is diagonal
in Kumjian's
sense.

\end{proof}

%
%
\bibliographystyle{amsplain}

\providecommand{\bysame}{\leavevmode\hbox to3em{\hrulefill}\thinspace}

\end{document}